\input amstex
\documentstyle{amsppt}
\topmatter
\title
Lattice Representations of Heisenberg Groups
\endtitle
\rightheadtext{Lattice Representations of Heisenberg Groups}
\author   Jae-Hyun Yang
\endauthor
\magnification =\magstep 1 \pagewidth{12.5cm} \pageheight{18.20cm}
\baselineskip =7mm
%Use \endgraf to indicate new paragraph
%\thanks will become a 1st page footnote
\thanks{This work was partially supported by the Max-Planck-Institut f\"{u}r
Mathematik and TGRC-KOSEF.\endgraf
}
\endthanks
\endtopmatter
\document
\NoBlackBoxes

\define\a{\alpha}
\define\be{\beta}

\define\G{\Gamma}

\define\vth{\vartheta}

\define\la{\lambda}

\define\s{\sigma}

\define\lrt{\longrightarrow}

\define\k{\kappa}
\define\ba{\backslash}

\define\M{{\Cal M}}

\define\BZ{\Bbb Z}
\define\BC{\Bbb C}
\define\BR{\Bbb R}

\define\N{N\in {\Bbb Z}^{(h,g)}}

\define\ex{\par\smallpagebreak\noindent}
\define\Box{$\square$}

\define\sd{\,\,{\vartriangleright}\kern -1.0ex{<}\,}
\define\sc{\bf}

\vskip 0.5cm
\head 1 \  Introduction %bold;
\endhead       %don't type final punctuation
\vskip 0.4cm
\ \ \
For any positive integers $g$ and $h$, we consider the Heisenberg group
$$H_{\BR}^{(g,h)}:=\left\{\,(\la,\mu,\kappa)\,\vert\ \la,\mu\in \BR^{(h,g)},\
\kappa\in \BR^{(h,h)},\ \kappa+\mu\,^t\!\la\ \text{symmetric}\ \right\}.$$
Recall that the multiplication law is
$$(\la,\mu,\kappa)\circ (\la',\mu',\kappa'):=(\la+\la',\mu+\mu',\kappa+
\kappa'+\la\,^t\!\mu'-\mu\,^t\!\la'). $$
Here $\BR^{(h,g)}\,(\,\text{resp.}\ \BR^{(h,h)})$ denotes the set
of all $h\times g\,(\,\text{resp.}\ h\times h)$ real matrices.
\par\medpagebreak\indent
The Heisenberg group $H_{\BR}^{(g,h)}$ is embedded to the symplectic group
$Sp(g+h,\BR)$ via the mapping
$$H_{\BR}^{(g,h)}\ni (\la,\mu,\kappa)\longmapsto
\pmatrix E_g & 0 & 0 & ^t\mu \\
\la & E_h & \mu & \kappa \\
0 & 0 & E_g & -^t\la \\
0 & 0 & 0 & E_h \endpmatrix \in Sp(g+h,\BR).$$
This Heisenberg group is a 2-step nilpotent Lie group and is important in
the study of toroidal compactifications of Siegel moduli spaces. In fact,
$H_{\BR}^{(g,h)}$ is obtained as the unipotent radical of the parabolic
subgroup of $Sp(g+h,\BR)$ associated with the rational boundary component
$F_g$\,(\,cf.\,[F-C]\,p.\,123 or [N] p.\,21\,). For the motivation of the study of
this Heisenberg group
we refer to [Y4]-[Y8] and [Z]. We refer to [Y1]-[Y3] for more results on
$H_{\BR}^{(g,h)}$. \par\medpagebreak
\indent
In [C], P. Cartier stated without proof that for $h=1,$ the lattice
representation of $H_{\BR}^{(g,1)}$ associated to the lattice $L$ is
unitarily equivalent to the direct sum of $[L^*:L]^{\frac 12}$ copies of
the Schr{\"{o}}dinger representation of $H_{\BR}^{(g,1)},$ where $L^*$ is
the dual lattice of $L$ with respect to a certain nondegenerate alternating
bilinear form. R. Berndt proved that the above fact for the case
$h=1$ in his lecture notes [B].
In this paper, we give a complete proof of Cartier's theorem
for $H_{\BR}^{(g,h)}.$
\vskip0.2cm
\noindent
{\bf Main\ Theorem.}\ \ Let $\M$ be a positive definite, symmetric
half-integral matrix of degree $h$ and $L$ be a self-dual lattice in
$\BC^{(h,g)}.$ Then the lattice representation $\pi_{\M}$
of $H_{\BR}^{(g,h)}$
associated with $L$ and $\M$ is unitarily equivalent to the direct sum
of $\left(\,\text{det}\,2\M\,\right)^g$ copies of the Schr{\"{o}}dinger
representation of $H_{\BR}^{(g,h)}.$ For more details, we refer to
Section 3.
\vskip 0.2cm
\indent
The paper is organized as follows. In Section 2, we review the
Schr{\"{o}}dinger representations of the Heisenberg group $H_{\BR}^{(g,h)}$.
In Section 3, we prove the main theorem. In the final section, we provide a relation
between lattice representations and theta functions.\par\medpagebreak
\noindent
{\smc Acknowledgement.} This work was in part done during my stay at the
Max-Planck-Institut f{\"u}r Mathematik in Bonn. I am very grateful to the
institute for hospitality and financial support. I also would like
to give my hearty thanks to the Department of Mathematics at Harvard University for its
hospitality during my short stay in Cambridge.
\par\medpagebreak
\noindent
{\smc Notations:}
\ \ We denote by $\BZ,\,\BR$ and $\BC$ the ring of integers, the field of real
numbers, and the field of complex numbers respectively. The symbol $\BC_1^{\times}$
denotes the multiplicative group consisting of all complex numbers $z$
with $\vert z\vert =1$, and the symbol
$Sp(g,\BR)$ the symplectic group of degree $g$, $H_g$
the Siegel upper half plane of degree $g$.
The symbol ``:='' means that the expression on the
right hand side is the definition of that on the left. We denote by $\BZ^+$
the set of all positive integers, by $F^{(k,l)}$ the set of
all $k\times l$ matrices with entries in a commutative ring $F$.
For any $M\in F^{(k,l)},\ ^t\!M$ denotes the transpose matrix of $M$.
For $A\in F^{(k,k)},\ \sigma(A)$ denotes the trace of $A$. For
$A\in F^{(k,l)}$ and $B\in F^{(k,k)},$ we set $B[A]={^t\!A}BA$.
We denote
the identity matrix of degree $k$ by $E_k$. For a positive integer $n,\ \text{Symm}\,(n,K)$
denotes the vector space consisting of all symmetric $n\times n$ matrices
with entries in a field $K.$ \ex
\vskip 0.5cm
\head  2 \  Schr\"{o}dinger Representations  %bold;
\endhead                              %don't type final punctuation
\vskip 0.5cm
\ \ \
First of all, we observe that
$H_{\BR}^{(g,h)}$ is a 2-step nilpotent Lie group. It is easy to see that
the inverse of an element $(\la,\mu,\kappa)\in H_{\BR}^{(g,h)}$ is given by
$$(\la,\mu,\kappa)^{-1}=(-\la,-\mu,-\kappa+\la\,^t\!\mu-\mu\,^t\!\la).$$
Now we set
$$[\la,\mu,\kappa]:=(0,\mu,\kappa)\circ (\la,0,0)=(\la,\mu,\kappa-\mu\,^t\!
\la).\tag2.1$$
Then $H_{\BR}^{(g,h)}$ may be regarded as a group
equipped with the following multiplication
$$[\la,\mu,\kappa]\diamond [\la_0,\mu_0,\kappa_0]:=[\la+\la_0,\mu+\mu_0,
\kappa+\kappa_0+\la\,^t\!\mu_0+\mu_0\,^t\!\la].\tag2.2$$
The inverse of $[\la,\mu,\kappa]\in H_{\BR}^{(g,h)}$ is given by
$$[\la,\mu,\kappa]^{-1}=[-\la,-\mu,-\kappa+\la\,^t\!\mu+\mu\,^t\!\la].$$
We set
$$K:=\left\{\,[0,\mu,\kappa]\in H_{\BR}^{(g,h)}\,\Big|\ \mu\in \BR^{(h,g)},\
\kappa=\,^t\!\kappa\in \BR^{(h,h)}\ \right\}.\tag2.3$$
Then $K$ is a commutative normal subgroup of $H_{\BR}^{(g,h)}$. Let
${\hat {K}}$ be the Pontrajagin dual of $K$, i.e., the commutative group
consisting of all unitary characters of $K$. Then ${\hat {K}}$ is
isomorphic to the additive group $\BR^{(h,g)}\times \text{Symm}\,(h,\BR)$ via
$$<a,{\hat {a}}>:=e^{2\pi i\s({\hat {\mu}}\,^t\!\mu+{\hat {\kappa}}\kappa)},
\ \ \ a=[0,\mu,\kappa]\in K,\ {\hat {a}}=({\hat {\mu}},{\hat {\kappa}})
\in {\hat {K}}.\tag2.4$$
We put
$$S:=\left\{\,[\la,0,0]\in H_{\BR}^{(g,h)}\,\Big|\ \la\in \BR^{(h,g)}\,
\right\}\cong \BR^{(h,g)}.\tag2.5$$
Then $S$ acts on $K$ as follows:
$$\alpha_{\la}([0,\mu,\kappa]):=[0,\mu,\kappa+\la\,^t\!\mu+\mu\,^t\!\la],
\ \ \ [\la,0,0]\in S.\tag2.6$$
It is easy to see that the Heisenberg group $\left( H_{\BR}^{(g,h)},
\diamond\right)$ is isomorphic to the semi-direct product $S\ltimes K$
of $S$ and $K$ whose multiplication is given by
$$(\la,a)\cdot (\la_0,a_0):=(\la+\la_0,a+\alpha_{\la}(a_0)),\ \ \la,\la_0\in S,\
a,a_0\in K.$$
On the other hand, $S$ acts on ${\hat {K}}$ by
$$\alpha_{\la}^{*}({\hat {a}}):=({\hat {\mu}}+2{\hat {\kappa}}\la,
{\hat {\kappa}}),\ \ [\la,0,0]\in S,\ \ a=({\hat {\mu}},{\hat {\kappa}})\in
{\hat {K}}.\tag2.7$$
Then, we have the relation $<\alpha_{\la}(a),{\hat {a}}>=<a,\alpha_{\la}^{*}
({\hat {a}})>$ for all $a\in K$ and ${\hat {a}}\in {\hat {K}}.$
\par\smallpagebreak
\indent
We have two types of $S$-orbits in ${\hat {K}}.$\par\smallpagebreak
\noindent
{\smc Type}
{\smc I.}\ \ Let ${\hat {\kappa}}\in \text{Symm}\,(h,\BR)$ with ${\hat {\kappa}}
\neq 0.$ The $S$-orbit of ${\hat {a}}({\hat {\kappa}}):=(0,{\hat {\kappa}})
\in {\hat {K}}$ is given by
$${\hat {\Cal O}}_{\hat {\kappa}}:=\left\{\,
(2{\hat {\kappa}}\la,{\hat {\kappa}})
\in {\hat {K}}\ \Big|\ \la\in \BR^{(h,g)}\,\right\}\,\cong\, \BR^{(h,g)}.
\tag2.8$$
{\smc Type} {\smc II.}\ \ Let ${\hat {y}}\in \BR^{(h,g)}.$ The $S$-orbit
${\hat {\Cal O}}_{\hat {y}}$ of ${\hat {a}}({\hat {y}}):=({\hat {y}},0)$
is given by
$${\hat {\Cal O}}_{\hat {y}}:=\left\{\,({\hat {y}},0)\,\right\}={\hat {a}}
({\hat {y}}).\tag2.9$$
We have
$${\hat {K}}=\left(\bigcup_{{\hat {\kappa}}\in \text{Symm}(h,\BR)}
{\hat {\Cal O}}_{\hat {\kappa}}\right)\bigcup \left( \bigcup_{{\hat {y}}\in
\BR^{(h,g)}}{\hat {\Cal O}}_{\hat {y}}\right) $$
as a set. The stabilizer $S_{\hat {\kappa}}$ of $S$ at ${\hat {a}}({\hat
{\kappa}})=(0,{\hat {\kappa}})$ is given by
$$S_{\hat {\kappa}}=\{0\}.\tag2.10$$
And the stabilizer $S_{\hat {y}}$ of $S$ at ${\hat {a}}({\hat {y}})=
({\hat {y}},0)$ is given by
$$S_{\hat {y}}=\left\{\,[\la,0,0]\,\Big|\ \la\in \BR^{(h,g)}\,\right\}=S
\,\cong\,\BR^{(h,g)}.\tag2.11$$
\indent
From now on, we set $G:=H_{\BR}^{(g,h)}$ for brevity.
It is known that $K$ is a closed,
commutative normal subgroup of $G$. Since $(\la,\mu,\kappa)=(0,\mu,
\kappa+\mu\,^t\!\la)\circ (\la,0,0)$ for $(\la,\mu,\kappa)\in G,$ the
homogeneous space $X:=K\backslash G$ can be identified with $\BR^{(h,g)}$
via
$$Kg=K\circ (\la,0,0)\longmapsto \la,\ \ \ g=(\la,\mu,\kappa)\in G.$$
We observe that $G$ acts on $X$ by
$$(Kg)\cdot g_0:=K\,(\la+\la_0,0,0)=\la+\la_0,\tag2.12$$
where $g=(\la,\mu,\kappa)\in G$ and $g_0=(\la_0,\mu_0,\kappa_0)\in G.$\ex
\indent
If $g=(\la,\mu,\kappa)\in G,$ we have
$$k_g=(0,\mu,\kappa+\mu\,^t\!\la),\ \ \ s_g=(\la,0,0)\tag2.13$$
in the Mackey decomposition of $g=k_g\circ s_g$\,(\,cf.\,[M]\,).
Thus if $g_0=(\la_0,\mu_0,\kappa_0)\in G,$ then we have
$$s_g\circ g_0=(\la,0,0)\circ (\la_0,\mu_0,\kappa_0)=(\la+\la_0,\mu_0,
\kappa_0+\la\,^t\!\mu_0)\tag2.14$$
and so
$$k_{s_g\circ g_0}=(0,\mu_0,\kappa_0+\mu_0\,^t\!\la_0+\la\,^t\!\mu_0+
\mu_0\,^t\!\la).\tag2.15$$
\indent
For a real symmetric matrix $c=\,^tc\in \BR^{(h,h)}$ with $c\neq 0$, we
consider the one-dimensional unitary representation $\sigma_c$ of $K$
defined by
$$\sigma_c\left((0,\mu,\kappa)\right):=e^{2\pi i\s(c\kappa)}\,I,\ \ \
(0,\mu,\kappa)\in K,\tag2.16$$
where $I$ denotes the identity mapping. Then the induced representation
$U(\s_c):=\text{ Ind}_K^G\,\s_c$ of $G$ induced from $\s_c$ is realized in the
Hilbert space ${\Cal H}_{\s_c}=L^2(X,d{\dot {g}},\BC)
\cong L^2\left(\BR^{(h,g)},
d\xi\right)$ as follows. If $g_0=(\la_0,\mu_0,\kappa_0)\in G$ and
$x=Kg\in X$ with $g=(\la,\mu,\kappa)\in G,$ we have
$$\left(U_{g_0}(\s_c)f\right)(x)=\s_c\left(k_{s_g\circ g_0}\right)\left(
f(xg_0)\right),\ \ f\in {\Cal H}_{\s_c}.\tag2.17$$
It follows from (2.15) that
$$\left( U_{g_0}(\s_c)f\right)(\la)=e^{2\pi i\s\{c(\kappa_0+\mu_0\,^t\!\la_0+
2\la\,^t\!\mu_0)\}}\,f(\la+\la_0).\tag2.18$$
Here, we identified $x=Kg$\,(resp.\,$xg_0=Kgg_0$) with $\la$\,(resp.\,
$\la+\la_0$). The induced representation $U(\s_c)$ is called the
{\it Schr{\" {o}}dinger\ representation} of $G$ associated with $\s_c.$
Thus $U(\s_c)$ is a monomial representation.\ex
\indent
Now, we denote by ${\Cal H}^{\s_c}$ the Hilbert space
consisting of all functions $\phi:G\lrt \BC$ which satisfy the following
conditions:\par \medpagebreak
\indent (1) $\phi(g)$ is measurable with respect to $dg,$\par\smallpagebreak
\indent (2) $\phi\left( (0,\mu,\kappa)\circ g)\right)=
e^{2\pi i\s(c\kappa)}\phi(g)$ \ \ for\ all $g\in G,$\par\smallpagebreak
\indent (3) $\parallel\phi\parallel^2:=\int_X\,\vert \phi(g)\vert^2\,
d{\dot {g}} < \infty,\ \ \ {\dot {g}}=Kg,$\par
\medpagebreak\noindent
where $dg$\,(resp.\,$d{\dot {g}}$) is a $G$-invariant measure on $G$
(resp.\,$X=K\backslash G$). The inner product $(\,,\,)$ on ${\Cal H}^{\s_c}$
is given by
$$(\phi_1,\phi_2):=\int_G\,\phi_1(g)\,{\overline {\phi_2(g)}}\,dg\ \ \ \
\text{for} \ \phi_1,\, \phi_2\in {\Cal H}^{\s_c}.$$
We observe that the mapping
$\Phi_c:{\Cal H}_{\s_c}\lrt {\Cal H}^{\s_c}$ defined by
$$\left( \Phi_c(f)\right)(g):=e^{2\pi i\s
\{c(\kappa+\mu\,^t\!\la)\}}\,f(\la),\ \ \ f\in {\Cal H}_{\s_c},\ g=(\la,\mu,
\kappa)\in G\tag 2.19$$
is an isomorphism of Hilbert spaces. The inverse $\Psi_c:{\Cal H}^{\s_c}
\lrt {\Cal H}_{\s_c}$ of $\Phi_c$ is given by
$$\left( \Psi_c(\phi)\right)(\la):=\phi((\la,0,0)),\ \ \
\phi\in {\Cal H}^{\s_c},\ \la\in \BR^{(h,g)}.\tag2.20$$
The Schr{\" {o}}dinger representation $U(\s_c)$ of $G$ on ${\Cal H}^{\s_c}$
is given by
$$\left( U_{g_0}(\s_c)\phi\right)(g)=
e^{2\pi i\s\{ c(\kappa_0+\mu_0\,^t\!\la_0+\la\,^t\!\mu_0-\la_0\,^t\!\mu)\}}\,
\phi\left( (\la_0,0,0)\circ g\right),\tag2.21$$
where $g_0=(\la_0,\mu_0,\kappa_0),\ g=(\la,\mu,\kappa)\in G$ and
$\phi\in {\Cal H}^{\s_c}.$ (2.21) can be expressed as follows.
$$\left(\,U_{g_0}(\s_c)\phi\,\right)(g)=e^{2\pi i\s\{c(\kappa_0+\kappa+
\mu_0\,^t\!\la_0+\mu\,^t\!\la+2\la\,^t\mu_0)\}}\,\phi((\la_0+\la,0,0)).
\tag2.22$$ \par\smallpagebreak
\noindent
{\sc Theorem\ 2.1.}\ Let $c$ be a positive symmetric
half-integral matrix of degree $h$. Then the Schr{\" {o}}dinger
representation $U(\s_c)$ of $G$ is irreducible.\par\medpagebreak
\noindent
{\it Proof.}\ The proof can be found in [Y1],\ theorem 3.
\hfill \Box\par
\bigpagebreak
\vskip 0.52cm
\head 3 \  Proof\ of\ the Main\ Theorem  %bold;
\endhead                             %don't type final punctuation
\vskip 0.24cm
Let $L:=\BZ^{(h,g)}\times \BZ^{(h,g)}$ be the lattice in the vector space
$V\cong \BC^{(h,g)}.$ Let $B$ be an alternating bilinear form on $V$ such
that $B(L,L)\subset \BZ,$ that is, $\BZ$-valued on $L\times L.$ The dual
$L_B^{*}$ of $L$ with respect to $B$ is defined by
$$L_B^{*}:=\left\{\,v\in V\,\vert\ B(v,l)\in \BZ\ \text{for\ all}\ l\in L\,
\right\}.$$
Then $L\subset L_B^{*}.$ If $B$ is nondegenerate, $L_B^{*}$ is also a
lattice in $V,$ called the {\it dual\ lattice} of $L$. In case $B$ is
nondegenerate, there exist a $\BZ$-basis $\{\,\xi_{11},\xi_{12},\cdots,
\xi_{hg},\eta_{11},\eta_{12},\cdots$,\par\noindent
$\eta_{hg}\,\}$ of $L$ and a set
$\{\,e_{11},e_{12},\cdots,e_{hg}\,\}$ of positive integers such
that
$e_{11}\vert e_{12},\,e_{12}\vert e_{13},$\par\noindent
$\cdots,e_{h,g-1}\vert e_{hg}$
for which
$$\pmatrix B(\xi_{ka},\xi_{lb})& B(\xi_{ka},\eta_{lb})\\
B(\eta_{ka},\xi_{lb}) & B(\eta_{ka},\eta_{lb}\endpmatrix=
\pmatrix 0&e\\ -e&0\endpmatrix,$$
where $1\leq k,l\leq h,\,1\leq a,b\leq g$ and $e:=\text{diag}\,(e_{11},e_{12},
\cdots,e_{hg})$ is the diagonal matrix of degree $hg$ with entries
$e_{11},e_{12},\cdots,e_{hg}.$ It is well known that $[L_B^{*}:L]=
(\,\text{det}\,e\,)^2=(e_{11}e_{12}\cdots e_{hg})^2$\,(cf.\,[I]\,p.\,72). The
number $\,\text{det}\,e\,$ is called the {\it Pfaffian} of $B.$\ex
\indent
Now, we consider the following subgroups of $G$:
$$\G_L:=\left\{\,(\la,\mu,\k)\in G\,\vert\ (\la,\mu)\in L,\ \k\in
\BR^{(h,h)}\,\right\}\tag3.1$$
and
$$\G_{L_B^{*}}:=\left\{\,(\la,\mu,\k)\in G\,\vert\ (\la,\mu)\in L_B^{*},\
\k\in \BR^{(h,h)}\,\right\}.\tag3.2$$
Then both $\G_L$ and $\G_{L_B^{*}}$ are normal subgroups of $G.$ We set
$${\Cal Z}_0:=\left\{\,(0,0,\k)\in G\,\vert\ \k=\,^t\k\in
\BZ^{(h,h)}\ \text{integral}\ \right\}.\tag3.3$$
It is easy to show that
$$\G_{L_B^{*}}=\left\{\,g\in G\,\vert\ g\gamma g^{-1}\gamma^{-1}\in
{\Cal Z}_0\ \text{for\ all}\ \gamma\in \G_L\,\right\}.$$
We define
$$Y_L:=\left\{\,\varphi\in {\text{ Hom}}\,(\G_L,\BC_1^{\times})\,\vert\
\varphi\ \text{is\ trivial\ on}\ {\Cal Z}_0\,\right\}\tag3.4$$
and
$$Y_{L,S}:=\left\{\,\varphi\in Y_L\,\vert\ \varphi(\k)=
e^{2\pi i\s(S\k)}\ \text{for
\ all}\ \k=\,^t\k\in \BR^{(h,h)}\,\right\}\tag3.5$$
for each symmetric real matrix $S$ of degree $h.$ We observe that, if $S$
is not half-integral, then $Y_L=\emptyset$ and so $Y_{L,S}=\emptyset.$ It
is clear that, if $S$ is symmetric half-integral, then $Y_{L,S}$ is
not empty.\ex
Thus we have
$$Y_L=\cup_{\M}\,Y_{L,\M},\tag3.6$$
where $\M$ runs through the set of all symmetric half-integral matrices of
degree $h$.\par\bigpagebreak
\noindent
{\bf Lemma\ 3.1.} Let $\M$ be a symmetric half-integral matrix of degree
$h$ with $\M\neq 0.$ Then any element $\varphi$ of $Y_{L,\M}$ is of the form
$\varphi_{\M,q}.$ Here $\varphi_{\M,q}$ is the character of $\G_L$ defined by
$$\varphi_{\M,q}((l,\k)):=e^{2\pi i\s(\M \k)}\cdot e^{\pi iq(l)}\ \
\ \text{for}\ (l,\k)\in \G_L,\tag3.7$$
where $q:L\lrt \BR/2\BZ\cong [0,2)$ is a function on $L$ satisfying the
following condition:
$$q(l_0+l_1)\equiv q(l_0)\,+\,
q(l_1)-2\s\{ \M(\la_0\,^t\!\mu_1-\mu_0\,^t\!\la_1)\}
\ \ (\,{\text{mod}}\ 2\,)\tag3.8$$
for all $l_0=(\la_0,\mu_0)\in L$ and $l_1=(\la_1,\mu_1)\in L.$
\par\medpagebreak
\noindent
{\it Proof.} (3.8) follows immediately from the fact that $\varphi_{\M,q}$ is
a character of $\G_L.$ It is obvious that any element of $Y_{L,\M}$ is of
the form $\varphi_{\M,q}.$ \hfill \Box\par\bigpagebreak
\noindent
{\bf Lemma\ 3.2.} An element of $Y_{L,0}$ is of the form $\varphi_{k,l}\,
(k,l\in \BR^{(h,g)}).$ Here $\varphi_{k,l}$ is the character
of $\G_L$ defined
by
$$\varphi_{k,l}(\gamma):=e^{2\pi i\s(k\,^t\!\la+l\,^t\!\mu)},\ \
\gamma=(\la,\mu,\k)\in \G_L.\tag3.9$$
\noindent
{\it Proof.} It is easy to prove and so we omit the proof.
\hfill \Box\par\bigpagebreak
\noindent
{\bf Lemma\ 3.3.} Let $\M$ be a nonsingular symmetric half-integral matrix
of degree $h$. Let $\varphi_{\M,q_1}$ and $\varphi_{\M,q_2}$ be the characters of
$\G_L$ defined by (3.7). The character $\varphi$ of $\G_L$ defined by
$\varphi:= \varphi_{\M,q_1}\cdot \varphi_{\M,q_2}^{-1}$
is an element of $Y_{L,0}.$\par\medpagebreak
\noindent
{\it Proof.} It follows from the existence of an element
$g=(\M^{-1}\la,\M^{-1}\mu,0)\in G$ with $(\la,\mu)\in V$ such that
$$\varphi_{\M,q_1}(\gamma)=\varphi_{\M,q_2}(g\gamma g^{-1})\ \
\text{for\ all}\
\gamma\in \G_L.$$
\hfill \Box\ex
\indent
For a unitary character $\varphi_{\M,q}$ of $\G_L$ defined by (3.7), we let
$$\pi_{\M,q}:={\text{ Ind}}_{\G_L}^G\,\varphi_{\M,q}\tag3.10$$
be the representation of $G$ induced from $\varphi_{\M,q}.$ Let
${\Cal H}_{\M,q}$ be the Hilbert space consisting of all measurable functions
$\phi:G\lrt \BC$ satisfying\ex
\indent (L1) $\phi(\gamma g)=\varphi_{\M,q}(\gamma)\,\phi(g)$ for all
$\gamma\in \G_L$ and $g\in G.$\par
\indent (L2)\ \ $\|\phi\|^2_{\M,q}:=\int_{\G_L\ba G}\,\vert \phi({\bar {g}})
\vert^2\,d{\bar {g}} <\infty,\ \ {\bar {g}}=\G_L g.$\ex
The induced representation $\pi_{\M,q}$ is realized in ${\Cal H}_{\M,q}$
as follows:
$$\biggl(\,\pi_{\M,q}(g_0)\phi\,\biggr)(g):=\phi(gg_0),\ \ g_0,g\in G,\
\phi\in {\Cal H}_{\M,q}.\tag3.11$$
The representation
$\pi_{\M,q}$ is called the {\it lattice\ representation} of $G$ associated
with the lattice $L$.\par\bigpagebreak
\noindent
{\bf Main\ Theorem.} Let $\M$ be a positive definite, symmetric half integral
matrix of degree $h$. Let $\varphi_{\M}$ be the character of $\G_L$ defined
by $\varphi_{\M}((\la,\mu,\k)):=e^{2\pi i\s(\M \k)}$ for all $(\la,\mu,\k)\in
\G_L.$ Then the lattice representation
$$\pi_{\M}:={\text{ Ind}}_{\G_L}^G\,\varphi_{\M}$$
induced from the character $\varphi_{\M}$ is unitarily equivalent to
the direct sum
$$\bigoplus\,U(\s_{\M}):=\bigoplus\,{\text{Ind}}_K^G\,\s_{\M}\ \ \
(\,(\,\hbox{$\text{det}\,2\M\,)^g$-copies}\,)$$
of the Schr{\"{o}}dinger representation ${\text{Ind}}_K^G\,\s_{\M}$.
\par\bigpagebreak
\noindent
{\it Proof.} We first recall that the induced representation $\pi_{\M}$ is
realized in the Hilbert space ${\Cal H}_{\M}$ consisting of all
measurable functions $\phi:G\lrt \BC$ satisfying the conditions
$$\phi((\la_0,\mu_0,\k_0)\circ g)=e^{2\pi i\s(\M \k_0)}\,\phi(g),\ \ \
(\la_0,\mu_0,\k_0)\in \G_L,\ g\in G \tag3.13$$
and
$$\|\phi\|^2_{\pi,\M}:=\int_{\G_L\ba G}\,\vert\phi({\bar g})\vert^2\,
d{\bar g} <\infty,\ \ {\bar g}=\G_L\circ g.\tag3.14$$
Now, we write
$$g_0=[\la_0,\mu_0,\k_0]\in \G_L\ \ and\ \ g=[\la,\mu,\k]\in G.$$
For $\phi\in {\Cal H}_{\M},$ we have
$$\phi(g_0\diamond g)=\phi([\la_0+\la,\mu_0+\mu,\k_0+\k+\la_0\,^t\!\mu+\mu\,
^t\!\la_0]).\tag3.15$$
On the other hand, we get
$$
\align
\phi(g_0\diamond g)&=\phi((\la_0,\mu_0,\k_0-\mu_0\,^t\!\la_0)\circ g)\\
&=e^{2\pi i\s\{ \M (\k_0-\mu_0\,^t\!\la_0)\}}\,\phi(g)\\
&=e^{2\pi i\s(\M \k_0)}\,\phi(g)\ \ \
(\,\text{because}\ \s(\M\mu_0\,^t\!\la_0)
\in \BZ\,)
\endalign$$
Thus, putting $\k^{\prime}:=\k_0+\la_0\,^t\!\mu+\mu\,^t\!\la_0,$ we get
$$\phi([\la_0+\la,\mu_0+\mu,\k+\k'])=e^{2\pi i\s(\M\k')}\cdot
e^{-4\pi i\s(\M \la_0\,^t\!\mu)}\,\phi([\la,\mu,\k]).\tag3.16$$
\indent
Putting $\la_0=\k_0=0$ in (3.16), we have
$$\phi([\la,\mu+\mu_0,\k])=\phi([\la,\mu,\k])\ \
\text{for\ all}\ \mu_0\in
\BZ^{(h,g)}\ \text{and}\ [\la,\mu,\k]\in G.\tag3.17$$
Therefore if we fix $\la$ and $\k,\ \phi$ is periodic in $\mu$ with respect to
the lattice $\BZ^{(h,g)}$ in $\BR^{(h,g)}.$ We note that
$$\phi([\la,\mu,\k])=\phi([0,0,\k]\diamond [\la,\mu,0])=e^{2\pi i\s(\M \k)}\,
\phi([\la,\mu,0])$$
for $[\la,\mu,\k]\in G.$ Hence, $\phi$ admits a Fourier expansion in $\mu:$
$$\phi([\la,\mu,\k])=e^{2\pi i\s(\M \k)}\sum_{\N}\,c_N(\la)\,
e^{2\pi i\s(N\,^t\!\mu)}.\tag3.18$$
If $\la_0\in \BZ^{(h,g)},$ then we have
$$\align
\phi([\la+\la_0,\mu,\k])&=e^{2\pi i\s(\M \k)}\sum_{\N}\,
c_N(\la+\la_0)\,e^{2\pi i\s(N\,^t\mu)}\\
&=e^{-4\pi i\s(\M \la_0\,^t\mu)}\,\phi([\la,\mu,\k])\ \ \
\ \ \ (\,\text{by}\ (3.16)\,)\\
&=e^{-4\pi i\s(\M \la_0\,^t\mu)}\,e^{2\pi i\s(\M\k)}\sum_{\N}
c_N(\la)\,e^{2\pi i\s(N\,^t\mu)},\\
&=e^{2\pi i\s(\M\k)}\,\sum_{\N}\,c_N(\la)\,
e^{2\pi i\s\{ (N-2\M\la_0)\,^t\mu\}}.\ \ \ (\,\text{by}\ (3.18)\,)
\endalign$$
So we get
$$\align
&  \sum_{\N}\,c_N(\la+\la_0)\,e^{2\pi i\s(N\,^t\mu)}\\
=&\sum_{\N}\,c_N(\la)\,e^{2\pi i\s\{ (N-2\M\la_0)\,^t\mu\}}\\
=&\sum_{\N}\,c_{N+2\M\la_0}(\la)\,e^{2\pi i\s(N\,^t\mu)}.
\endalign$$
Hence, we get
$$c_N(\la+\la_0)=c_{N+2\M\la_0}(\la)\ \
\text{for\ all}\ \la_0\in \BZ^{(h,g)}\ \text{and}\
\la\in \BR^{(h,g)}.\tag3.19$$
Consequently, it is enough to know only the coefficients $c_{\a}(\la)$
for the representatives $\a$ in $\BZ^{(h,g)}$ modulo $2\M$. It is obvious
that the number of all such $\a$'s is $(\text{det}\,2\M)^g.$
We denote by ${\Cal J}$ a complete system of such representatives
in $\BZ^{(h,g)}$ modulo $2\M.$\ex
Then, we have
$$\align
\ \phi([\la,\mu,\k])
=e^{2\pi i\s(\M\k)}\,\ &\biggl\{\, \sum_{\N}\,c_{\a+2\M N}(\la)\,
e^{2\pi i\s\{ (\a+2\M N)\,^t\mu\}}\\
&+\sum_{\N}\,c_{\be+2\M N}(\la)\,e^{2\pi i\s\{ (\be+2\M N)\,^t\mu\}}\\
&\cdot\\
&\cdot\\
&+\sum_{\N}\,c_{\gamma+2\M N}(\la)\,e^{2\pi i\{ (\gamma+2\M N)\,^t\mu\}}
\,\biggr\} ,
\endalign$$
where $\{\,\a,\be,\cdots,\gamma\,\}$ denotes the complete system
${\Cal J}.$\ex
\indent
For each $\a\in {\Cal J},$ we denote by ${\Cal H}_{\M,\a}$ the Hilbert
space consisting of Fourier expansions
$$e^{2\pi i\s(\M \k)}\,\sum_{\N}\,c_{\a+2\M N}(\la)\,
e^{2\pi i\s\{\,(\a+2\M N)\,^t\mu\}},\ \ \ (\la,\mu,\k)\in G,$$
where $c_N(\la)$ denotes the coefficients of the Fourier expansion (3.18) of
$\phi\in {\Cal H}_{\M}$ and $\phi$ runs over the set $\{\,\phi\in
\pi_{\M}\,\}$. It is easy to see that ${\Cal H}_{\M,\a}$ is invariant under
$\pi_{\M}.$ We denote the restriction of $\pi_{\M}$ to ${\Cal H}_{\M,\a}$
by $\pi_{\M,\a}.$ Then we have
$$\pi_{\M}=\bigoplus_{\a\in {\Cal J}}\,\pi_{\M,\a}.\tag3.20$$
Let $\phi_{\a}\in \pi_{\M,\a}.$ Then for $[\la,\mu,\k]\in G,$ we get
$$\phi_{\a}([\la,\mu,\k])=e^{2\pi i\s(\M \k)}\,\sum_{\N}\,
c_{\a+2\M N}(\la)\,e^{2\pi i\s\{ (\a+2\M N)\,^t\mu\}}.\tag3.21$$
We put
$$I_{\la}:=\overbrace{[0,1]\times [0,1]\times \cdots \times [0,1]}^{(h\times
g)\text{-times}}
\, \subset \left\{\,[\la,0,0]\,\vert\ \la\in \BR^{(h,g)}\,\right\}$$
and
$$I_{\mu}:=\overbrace{[0,1]\times [0,1]\times \cdots \times [0,1]}^{(h\times
g)\text{-times}}
\, \subset\, \left\{\,[0,\mu,0]\,\vert\ \mu\in \BR^{(h,g)}\right\}.$$
Then, we obtain
$$\int_{I_{\mu}}\,\phi_{\a}([\la,\mu,\k])\,e^{-2\pi i\s(\a\,^t\mu)}\,d\mu
=e^{2\pi i\s(\M \k)}\,c_{\a}(\la),\ \ \a\in {\Cal J}.\tag3.22$$
Since $\G_L\ba G\cong I_{\la}\times I_{\mu},$ we get
$$\align
\|\phi_{\a}\|^2_{\pi,\M,\a}:&=\|\phi_{\a}\|^2_{\pi,\M}=
\int_{\G_L\ba G}\,\vert \phi_{\a}({\bar g})\vert^2\,d{\bar g}\\
&=\int_{I_{\la}}\int_{I_{\mu}}\,\vert \phi_{\a}({\bar g})\vert^2\,d\la d\mu\\
&=\int_{I_{\la}\times I_{\mu}}\biggl|\sum_{\N}\,
c_{\a+2\M N}(\la)\,e^{2\pi i\s\{ (\a+2\M N)\,^t\mu\}}\,\biggr|^2\,
d\la d\mu\\
&=\int_{I_{\la}}\,\sum_{\N}\,\vert c_{\a+2\M N}(\la)\vert^2\, d\la\\
&=\int_{I_{\la}}\,\sum_{\N}\,\vert c_{\a}(\la+N)\vert^2\,d\la\ \ \
\ \ \ (\,\text{by}\ (3.19)\,)\\
&=\int_{\BR^{(h,g)}}\,\vert c_{\a}(\la)\vert^2\,d\la.
\endalign$$
Since $\phi_{\a}\in \pi_{\M,\a},\ \|\phi_{\a}\|_{\pi,\M,\a}<\infty$ and
so $c_{\a}(\la)\in L^2\left(\BR^{(h,g)},d\xi\right)$ for all
$\a\in {\Cal J}.$\ex
\indent
For each $\a\in {\Cal J},$ we define the mapping $\vth_{\M,\a}$ on
$L^2\left( \BR^{(h,g)},d\xi\right)$ by
$$(\vth_{\M,\a}f)([\la,\mu,\k]):=e^{2\pi i\s(\M \k)}\,\sum_{\N}\,
f(\la+N)\,e^{2\pi i\s\{ (\a+2\M N)\,^t\mu\}},\tag3.23$$
where $f\in L^2\left( \BR^{(h,g)},d\xi\right)$ and $[\la,\mu,\k]\in G.$
\par\bigpagebreak
\noindent
{\bf Lemma\ 3.4.} For each $\a\in {\Cal J},$ the image of
$L^2\left( \BR^{(h,g)},d\xi\right)$ under $\vth_{\M,\a}$ is contained in
${\Cal H}_{\M,\a}.$ Moreover, the mapping
$\vth_{\M,\a}$ is a one-to-one unitary
operator of $L^2\left(\BR^{(h,g)},d\xi\right)$ onto ${\Cal H}_{\M,\a}$
preserving the norms. In other words, the mapping
$$\vth_{\M,\a}:L^2\left( \BR^{(h,g)},d\xi\right)\lrt {\Cal H}_{\M,\a}$$
is an isometry.\par\medpagebreak
\noindent
{\it Proof.} We already showed that $\vth_{\M,\a}$ preserves the norms.
First, we observe that if $(\la_0,\mu_0,\k_0)\in \G_L$ and $g=[\la,\mu,\k]
\in G,$
$$\align
(\la_0,\mu_0,\k_0)\circ g&=[\la_0,\mu_0,\k_0+\mu_0\,^t\!\la_0]\diamond
[\la,\mu,\k]\\
&=[\la_0+\la,\mu_0+\mu,\k+\k_0+\mu_0\,^t\!\la_0+\la_0\,^t\mu+\mu\,^t\!\la_0].
\endalign$$
Thus we get
$$\align
& (\vth_{\M,\a}f)((\la_0,\mu_0,\k_0)\circ g)\\
&=e^{2\pi i\s\{ \M(\k+\k_0+\mu_0\,^t\!\la_0+\la_0\,^t\mu+\mu\,^t\!\la_0)\}}\,
\sum_{\N}\,f(\la+\la_0+N)\,e^{2\pi i\{ (\a+2\M N)\,^t(\mu_0+\mu)\}}\\
&=e^{2\pi i\s (\M \k_0)}\cdot e^{2\pi i\s (\M \k)}
\cdot e^{2\pi i\s(\a\,^t\!\mu_0)}\,
\sum_{\N}\,f(\la+N)\,e^{2\pi i\s\{ (\a+2\M N)\,^t\mu\}}\\
&=e^{2\pi i\s (\M \k_0)}\,(\vth_{\M,\a}f)(g).
\endalign$$
Here, in the above equalities we used the facts that
$2\s(\M N\,^t\mu_0)\in \BZ$ and $\a\,^t\mu_0\in \BZ.$ It is easy to show
that
$$\int_{\G_L\ba G}\,\vert \vth_{\M,\a}f({\bar g})\vert^2\,d{\bar g}=
\int_{\BR^{(h,g)}}\,\vert f(\la)\vert^2\,d\la=\|f\|^2_2 < \infty.$$
This completes the proof of Lemma 3.4.\ \ex
\indent
Finally, it is easy to show that for each $\a\in {\Cal J},$ the mapping
$\vth_{\M,\a}$ intertwines the Schr{\" {o}}dinger representation
$\left( U(\s_{\M}),L^2(\BR^{(h,g)},d\xi)\right)$ and the representation
$(\pi_{\M,\a},{\Cal H}_{\M,\a}).$ Therefore, by Lemma 3.4, for each
$\a\in {\Cal J},\ \pi_{\M,\a}$ is unitarily equivalent to $U(\s_{\M})$ and
so $\pi_{\M,\a}$ is an irreducible unitary representation of $G.$
According to (3.20), the induced representation $\pi_{\M}$ is unitarily
equivalent to
$$\bigoplus\,U(\s_{\M})\ \ \ (\,(\,\hbox{$\text{det}\,2\M)^g$-copies}\,).$$
This completes the proof of the Main Theorem. \hfill \Box

\vskip 0.5cm
\head 4 \ Relation of Lattice Representations to Theta Functions %bold
\endhead %don't type final punctuation
\vskip 0.2cm
\indent
In this section, we state the connection between
lattice representations and theta
functions. As before, we write $\,V=\BR^{(h,g)}\times \BR^{(h,g)}\cong
\BC^{(h,g)},\ L=\BZ^{(h,g)}\times \BZ^{(h,g)}$ and $\M$ is a positive
symmetric half-integral matrix of degree $h$. The function
$q_{\M}:L\lrt \BR / 2\BZ\,=\,[0,2)$ defined by
$$q_{\M}((\xi,\eta)):=\,2\,\s(\M \xi\,^t\eta),\ \ \ (\xi,\eta)\in L
\tag4.1$$
satisfies Condition (3.8). We let $\,\varphi_{\M,q_{\M}}:\G_L\lrt
\BC_1^{\times}$ be the character of $\G_L$ defined by
$$\varphi_{\M,q_{\M}}((l,\kappa))\,=\,
e^{2\pi i\,\s(\M\kappa)}\,e^{\pi i\,q_{\M}(l)}\,,\ \ \ (l,\kappa)\in
\G_L.$$
We denote by ${\Cal H}_{\M,q_{\M}}$ the Hilbert space consisting of
measurable functions $\phi:G\lrt \BC$ which satisfy Condition (4.2)
and Condition (4.3):\par\medpagebreak
\noindent
(4.2)\ \ \ $\phi((l,\kappa)\circ g)\,=\,\varphi_{\M,q_{\M}}((l,\kappa))\,
\phi(g)$ \ \ for all $(l,\kappa)\in \G_L$ and $g\in G.$\par\smallpagebreak
$$\int_{{\G_L}\ba G}\,\| \phi({\dot g})\|^2\,d{\dot g}\,< \,\infty,\ \ \
{\dot g}=\G_L\circ g.\tag4.3$$
Then the lattice representation
$$\pi_{\M,q_{\M}}:=\text{ Ind}_{\G_L}^G\,\varphi_{\M,q_{\M}}$$
of $G$ induced from the character $\varphi_{\M,q_{\M}}$ is realized in
${\Cal H}_{\M,q_{\M}}$ as
$$\left(\,\pi_{\M,q_{\M}}(g_0)\,\phi\,\right)(g)\,=\,\phi(gg_0),\ \ \
g_0,g\in G,\ \phi\in {\Cal H}_{\M,q_{\M}}.$$
Let ${\bold H}_{\M,q_{\M}}$ be the vector space consisting of measurable functions
$\,F:V\lrt \BC\,$ satisfying Conditions (4.4) and (4.5).
$$F(\la+\xi,\mu+\eta)\,=\,e^{2\pi i\,\s\{ \M(\xi\,^t\eta+\la\,^t\eta-\mu\,
^t\xi)\} }\,F(\la,\mu)\tag4.4$$
for all $(\la,\mu)\in V$ and $(\xi,\eta)\in L.$
$$\int_{L\ba V}\,\| F({\dot v})\|^2\,d{\dot v}\,=\,
\int_{I_{\la}\times I_{\mu}}\,\| F(\la,\mu)\|^2\,d\la d\mu\,<\,\infty.
\tag4.5$$
\indent
Given $\,\phi\in {\Cal H}_{\M,q_{\M}}$ and a fixed element $\Omega\in H_g,$
we put\par\medpagebreak\noindent
(4.6)\ \ \ \ \ $E_{\phi}(\la,\mu):=\,\phi((\la,\mu,0)),\ \ \ \la,\mu\in
\BR^{(h,g)},$\par\smallpagebreak\noindent
(4.7)\ \ \ \ \ $F_{\phi}(\la,\mu):=\,\phi([\la,\mu,0]),\ \ \ \la,\mu\in
\BR^{(h,g)},$\par\smallpagebreak\noindent
(4.8)\ \ \ \ \ $F_{\Omega,\phi}(\la,\mu):=\,e^{-2\pi i\,\s(\M \la\Omega\,^t\!\la)}
\,F_{\phi}(\la,\mu),\ \ \ \la,\mu\in \BR^{(h,g)}.$\par\medpagebreak
\indent
In addition, we put for $W=\la\Omega+\mu\in \BC^{(h,g)},$
$$\vth_{\Omega,\phi}(W):=\,\vth_{\Omega,\phi}(\la\Omega+\mu):=\,F_{\Omega,\phi}
(\la,\mu).\tag4.9$$
We observe that $\,E_{\phi},\,F_{\phi}$ and $F_{\Omega,\phi}\,$ are functions
defined on $V$ and $\,\vth_{\Omega,\phi}\,$ is a function defined on
$\BC^{(h,g)}.$\par\bigpagebreak
\noindent
{\bf Proposition\ 4.1.} If $\,\phi\in {\Cal H}_{\M,q_{\M}},\ (\xi,\eta)\in L$
and $(\la,\mu)\in V,$ then we have the formulas
$$
 E_{\phi}(\la+\xi,\mu+\eta)\,=\,
e^{2\pi i\,\s\{ \M(\xi\,^t\eta+\la\,^t\eta-\mu\,^t\xi)\} }\,E_{\phi}(\la,\mu).
 \tag4.10$$
$$ F_{\phi}(\la+\xi,\mu+\eta)\,=\,
e^{ -4\pi i\,\s(\M\xi\,^t\mu) }\,F_{\phi}(\la,\mu).\tag4.11\hskip 10cm$$
$$ F_{\Omega,\phi}(\la+\xi,\mu+\eta)\,=\,
e^{ -2\pi i\,\s\{ \M(\xi\Omega\,^t\xi+2\la\Omega\,^t\xi+2\mu\,^t\xi)\} }\,
F_{\Omega,\phi}(\la,\mu).\tag4.12$$
If $\,W=\la\Omega+\eta\in \BC^{(h,g)},$ then we have
$$\vth_{\Omega,\phi}(W+\xi\Omega+\eta)\,=\,
e^{ -2\pi i\,\s\{ \M(\xi\Omega\,^t\xi+2W\,^t\xi)\} }\,\vth_{\Omega,\phi}
(W).\tag4.13$$
Moreover, $\,F_{\phi}\,$ is an element of ${\bold H}_{\M,q_{\M}}.$
\par\medpagebreak
\noindent
{\it Proof.} We note that
$$(\la+\xi,\mu+\eta,0)=(\xi,\eta,-\xi\,^t\mu+\eta\,^t\!\la)\circ
(\la,\mu,0).$$
Thus we have
$$\align
E_{\phi}(\la+\xi,\mu+\eta)&=\phi((\la+\xi,\mu+\eta,0))\\
&=\phi((\xi,\eta,\,-\xi\,^t\mu+\eta\,^t\!\la)\circ (\la,\mu,0))\\
&=\,e^{2\pi i\s\{ \M(\xi\,^t\eta+\la\,^t\eta-\mu\,^t\xi)\} }\,
\phi((\la,\mu,0))\\
&=\,e^{ 2\pi i\s\{ \M(\xi\,^t\eta+\la\,^t\eta-\mu\,^t\xi)\} }\,E_{\phi}
(\la,\mu).
\endalign$$
This proves Formula (4.10). We observe that
$$[\la+\xi,\mu+\eta,0]=(\xi,\eta,\,-\xi\,^t\mu-\mu\,^t\xi-\eta\,^t\xi)\circ
[\la,\mu,0].$$
Thus we have
$$\align
F_{\phi}(\la+\xi,\mu+\eta)&=\phi([\la+\xi,\mu+\eta,0])\\
&=\,e^{ -2\pi i\,\s\{ \M(\xi\,^t\mu+\mu\,^t\xi+\eta\,^t\xi)\} }\\
& \ \ \ \ \quad\times e^{2\pi i\,\s(\M\xi\,^t\eta)}\,\phi([\la,\mu,0])\\
&=\,e^{ -4\pi i\,\s(\M\xi\,^t\mu) }\,\phi([\la,\mu,0])\\
&=\,e^{ -4\pi i\,\s(\M\xi\,^t\mu) }\,F_{\phi}(\la,\mu).
\endalign$$
This proves Formula (4.11). According to (4.11), we have
$$\align
F_{\Omega,\phi}(\la+\xi,\mu+\eta)&=
\,e^{-2\pi i\,\s\{ \M (\la+\xi)\Omega\,^t(\la+\xi)\} }\,F_{\phi}(\la+\xi,\mu+
\eta)\\
&=\, e^{ -2\pi i\,\s\{ \M (\la+\xi)\Omega\,^t(\la+\xi)\} }\\
&  \ \ \ \ \quad\times e^{ -4\pi i\,\s(\M\xi\,^t\mu) }\,F_{\phi}(\la,\mu)\\
&=\,e^{ -2\pi i\,\s\{ \M (\xi\Omega\,^t\xi+2\la\Omega\,^t\xi+2\mu\,^t\xi)\} }
\\
& \ \ \ \ \quad\times e^{ -2\pi i\,\s(\M\la\Omega\,^t\!\la) }\,F_{\phi}(\la,\mu)\\
&=\,e^{ -2\pi i\,\s\{ \M(\xi\Omega\,^t\xi+2\la\Omega\,^t\xi+2\mu\,^t\xi)\} }
\, F_{\Omega,\phi}(\la,\mu).
\endalign$$
This proves Formula (4.12). Formula (4.13) follows immediately from
Formula (4.12). Indeed, if $\,W=\la\Omega+\mu\,$ with $\la,\mu\in \BR^{(h,g)},$ we have
$$\align
\vth_{\Omega,\phi}(W+\xi\Omega+\eta)&=F_{\Omega,\phi}(\la+\xi,\mu+\eta)\\
&=\,e^{ -2\pi i\,\s\{ \M(\xi\Omega\,^t\xi+2(\la\Omega+\mu)\,^t\xi)\} }\,
F_{\Omega,\phi}(\la,\mu)\\
&=\,e^{ -2\pi i\,\s\{ \M(\xi\Omega\,^t\xi+2W\,^t\xi)\} }\,
\vth_{\Omega,\phi}(W).
\endalign$$
\hfill \Box\ex
\noindent
{\smc Remark\ 4.2.} The function $\,\vth_{\Omega,\phi}(W)\,$ is a theta
function of level $2\M$ with respect to $\Omega$ if $\,\vth_{\Omega,\phi}\,$
is holomorphic. For any $\,\phi\in {\Cal H}_{\M,q_{\M}},$ the function
$\,\vth_{\Omega,\phi}\,$ satisfies the
well known transformation law of a theta
function. In this sense, the lattice representation
$\,(\,\pi_{\M,q_{\M}},\,{\Cal H}_{\M,q_{\M}}\,)$ is closely related to
theta functions.

\vskip 1.5cm \noindent Max-Planck Institut f{\" u}r
Mathematik\par\noindent Gottfried-Claren-Strasse 26 \par\noindent
D-53225 Bonn\par\noindent Germany \vskip 1.0cm {\it The\ present\
address\ is} \vskip 0.5cm \noindent Department of
Mathematics\par\noindent Inha University \par\noindent Inchon
402-751 \par\noindent Republic of Korea\par\medpagebreak \vskip
0.1cm \noindent {\smc E-mail\,:\ jhyang\@inha.ac.kr}

\end{document}